\documentclass[10pt,reqno]{amsart}
  \usepackage{geometry}
\newtheorem*{rep@theorem}{\rep@title}
\newcommand{\newreptheorem}[2]{%
\newenvironment{rep#1}[1]{%
 \def\rep@title{#2 \ref{##1}}%
 \begin{rep@theorem}}%
 {\end{rep@theorem}}}
\makeatother

 \usepackage{amscd,amsmath,amssymb,amsthm,amsfonts}
\usepackage[alphabetic]{amsrefs}
\usepackage{mathrsfs}
\usepackage[shortlabels]{enumitem}
\usepackage[all]{xy}

\usepackage{thmtools}

\usepackage{xcolor}
\usepackage[hypertexnames=true,colorlinks = true, allcolors=blue,linktoc = all, pdffitwindow = false, urlbordercolor = white]{hyperref}%

\usepackage{tikz}
\usetikzlibrary{decorations.pathreplacing}
 
\usetikzlibrary{knots}
\usetikzlibrary{decorations.markings}
\usepackage{hyperref}
\usepackage{graphicx} 

 \newcommand{\IN}[0]{\mathbb{N}}

 \newcommand{\IZ}[0]{\mathbb{Z}}

 \newcommand{\CF}[0]{\mathcal{F}}

 \newcommand{\CL}[0]{\mathcal{L}}

 \newcommand{\CT}[0]{\mathcal{T}}

\newcommand {\stab}{{\rm Stab}}

\newreptheorem{theorem}{Theorem}
\newtheorem{theorem}{Theorem}
\newtheorem*{theorem*}{Theorem}
\newtheorem*{question*}{Question}
\newtheorem*{proposition*}{Proposition}

\newtheorem*{lemma*}{Lemma}

\newtheorem{question}{Question}

\newtheorem{remark}[theorem]{Remark}

\newtheorem{convention}{Convention}[section]

\numberwithin{equation}{section}

\begin{document}
\title[Remarks on some maximal subgroups of $F$ and on the $\vec{F}$-index]
{Remarks on some maximal subgroups of $F$ and on the $\vec{F}$-index of knots}
\author{Valeriano Aiello} 
\address{Valeriano Aiello,
Dipartimento di Matematica, Universit\`a di Roma Sapienza, P.le Aldo Moro
5, 00185 Roma, Italy, \url{https://github.com/valerianoaiello}}
\email{valerianoaiello@gmail.com}

\begin{abstract}
We demonstrate that three maximal subgroups of infinite index in the rectangular subgroup \( K_{(2,2)} \) of the Thompson group \( F \), each containing Jones's \( 3 \)-colorable subgroup \( \mathcal{F} \), can be characterized as stabilizer subgroups. 
Additionally, we show that the \( \vec{F} \)-index, an elementary knot invariant introduced thanks to Jones's construction of knots from Thompson groups, 
  may increase at most by $4$  after changing the orientation of a knot.
  
  \bigskip
  
  \end{abstract}
\maketitle

\section*{Introduction}

Thompson’s group $F$ was introduced by R. Thompson in the 1960s, alongside its two sibling groups, $T$ and $V$. The group $F$ consists of piecewise linear homeomorphisms of $[0,1]$ with derivatives that are powers of $2$ and a finite number of points of non-differentiability at dyadic rationals. The group $T$ extends this action to the circle $S^1$, while $V$ acts on the Cantor set. Thompson’s group $F$ is one of the most intriguing countable discrete groups, yet it remains enigmatic, as its analytical properties have posed challenges to experts for decades. 
In particular, the question of its amenability is a notorious open problem.

In the 1980s, Brown introduced in \cite{Brown} a family of generalizations of $F$, known as the Brown-Thompson groups $F_p$ for $p \geq 2$ (see also \cite{Wladis} for a nice presentation of these groups). Elements of $F_p$ can still be viewed as piecewise linear homeomorphisms of $[0,1]$, but their derivatives are now powers of $p$, and the points of non-differentiability take the form $a/p^k$, where $a, k \in \mathbb{Z}$. In particular, for $p=2$, one recovers Thompson’s group, i.e., $F_2 = F$. 

Thompson groups have been widely studied across various mathematical fields. In recent years, they have attracted interest in cryptography \cite{ShCrypo, MatCrypto}, quantum probability \cite{Kostler, KKW, Kri, A3}, and knot theory \cite{Jo18, A, Kodama3, KodamaP, Krushkal}. In the latter context, Jones introduced a method to construct knots and links from elements of $F$ \cite{Jo14}. The oriented knots and links arise from a subgroup of $F$ known as the oriented subgroup $\vec{F}$. 
Originally Jones had proved that all oriented links arise from elements of $\vec{F}$ up to disjoint union with unlinks.
Later in \cite{A} it was proved that links could be exactly reproduced.
Consequently, $F$ and $\vec{F}$ can be viewed as alternatives to the braid groups. Notably, $F$ contains a monoid of positive elements, denoted $F_+$, while $\vec{F}$ contains a corresponding monoid $\vec{F}_+$, playing a role analogous to that of positive braids in the unoriented and oriented constructions. In the framework of Jones's construction of knots, these monoids were investigated in \cite{AB1, AB2}.

The oriented subgroup $\vec{F}$ has led to significant discoveries in subgroup structure of $F$, notably providing the first explicit example of a maximal subgroup of infinite index in $F$ that is not a stabilizer of a point \cite{GS, GS2}. Further developments based on Jones’s technology have produced a variety of interesting subgroups \cite{TV, TV2, TV3, Bro1, Cdl}, many of which are relevant to the study of maximal subgroups - a topic that has received growing attention in recent years (see, e.g., \cite{BBQS, SavBak, Gri}).

Since Jones’s seminal work on this topic \cite{Jo14}, a number of articles have expanded upon his ideas. For a more comprehensive overview, we refer the reader to the survey articles \cite{Bro0, A2}.

We end this introduction by presenting the main results of this article. 
We refer to Sections \ref{sec:1} and \ref{section_main_res} for the notation.

Among the aforementioned maximal subgroups, three subgroups $M_0$, $M_1$, and $M_2$ of $F$ were defined in \cite{TV2} by using Jones's $3$-colorable subgroup $\CF$. All three sit inside Bleak and Wassink's rectangular subgroup $K_{(2,2)}$. The latter was introduced, along with a family of subgroups of $F$, and studied in \cite{BW}. $K_{(2,2)}$ consists of the homeomorphisms $g$ of $F$ such that $\log_2 g'(0)$ and $\log_2 g'(1)$ are even integers.

The subgroups $M_0$, $M_1$, and $M_2$ were originally defined by means of a finite set of generators. The goal of this paper is to characterize them in terms of stabilizers of subsets of dyadic integers under the aforementioned natural action on $[0,1]$, namely the theorem below (see Section \ref{section_main_res} for the definitions of $M_0$, $M_1$, $M_2$, and of the subsets $S_0$, $S_1$, $S_2$).

\begin{theorem} 
The subgroups $M_0$, $M_1$, and $M_2$ can be characterized as
\[
M_0 = \stab(S_2), \quad 
M_1 = \stab(S_0) \cap K_{(2,2)}, \quad
M_2 = \stab(S_1) \cap K_{(2,2)},
\]
where $S_0$, $S_1$, $S_2$ are subsets of the dyadic rationals in $[0,1]$. 
\end{theorem}

The second main result has to do with the knot-theoretical applications obtained thanks to yet another maximal subgroup, the oriented subgroup $\vec{F}$. As with the braid groups, infinitely many elements produce the same link. One may consider the $\vec{F}$-index of an oriented knot $K$, namely the smallest number of leaves required by each binary tree in a binary tree diagram $(T_+,T_-)$ of $\vec{F}$ such that $K$ is realized as $\vec{\CL}(T_+,T_-)$.
To any oriented knot $K$ we can associate its reversed knot $K^r$, namely the knot which is the same topologically but with opposite orientation. 

\begin{theorem} 
Let $K$ be an oriented knot and $K^r$ be the corresponding reversed knot. 
Then the $\vec{F}$-indices of $K$ and $K^r$ differ by at most $4$.
\end{theorem}

\section{Preliminaries}\label{sec:1}

This section is devoted to revisiting the definitions of Thompson's group $F$, the Brown-Thompson groups $F_p$, the oriented subgroup $\vec{F}$ and Jones's construction of knots from elements of Thompson groups, and the $3$-colorable subgroup $\CF$. Readers interested in further details can refer to \cite{CFP, B, Brown, Wladis} for $F$ and $F_p$, \cite{GS} for $\vec{F}$,  \cite{Jo18, A2} for the construction of knots, and \cite{Ren,  TV2, TV3} for $\CF$.

Thompson's group $F$ consists of all piecewise linear homeomorphisms of the unit interval $[0,1]$ that are differentiable everywhere except at a finite set of dyadic rational numbers (i.e., numbers of the form $a/2^k$, where $a, k \in \mathbb{Z}$). Moreover, where the derivative is well-defined, it is a power of $2$. We adopt the conventional notation: for $f, g \in F$, the composition is given by $f \cdot g(t) = g(f(t))$.

Thompson's group $F$ admits the following infinite presentation:
\[
F = \langle x_0, x_1, x_2, \ldots \; | \; x_n x_k = x_k x_{n+1}, \quad \forall \; k<n \rangle.
\]
Notably, $x_0$ and $x_1$ are sufficient to generate $F$.

The elements $x_0$, $x_1$, $x_2, \ldots$ generate a monoid, denoted by $F_+$, whose elements are referred to as positive. It is worth mentioning that Thompson's monoid $F_+$ and Thompson-like monoids have been of interest in non-commutative probability; see, for example, \cite{KKW, CDGR, CDR}.

\begin{figure}
\phantom{This text will be invisible} 
\[
\begin{tikzpicture}[x=.35cm, y=.35cm,
    every edge/.style={
        draw,
      postaction={decorate,
                    decoration={markings}
                   }
        }
]

\node at (-1.5,0) {$\scalebox{1}{$x_0=$}$};
\node at (-1.25,-3) {\;};

\draw[thick] (0,0) -- (2,2)--(4,0)--(2,-2)--(0,0);
 \draw[thick] (1,1) -- (2,0)--(3,-1);

 \draw[thick] (2,2)--(2,2.5);

 \draw[thick] (2,-2)--(2,-2.5);

 \draw[thick, gray, dashed] (-.25,0)--(4.5,0);

\end{tikzpicture}\qquad
\;\;
\begin{tikzpicture}[x=.35cm, y=.35cm,
    every edge/.style={
        draw,
      postaction={decorate,
                    decoration={markings}
                   }
        }
]

\node at (-3.5,0) {$\scalebox{1}{$x_1=$}$};
\node at (-1.25,-3.25) {\;};

\draw[thick] (2,2)--(1,3)--(-2,0)--(1,-3)--(2,-2);

\draw[thick] (0,0) -- (2,2)--(4,0)--(2,-2)--(0,0);
 \draw[thick] (1,1) -- (2,0)--(3,-1);

 \draw[thick] (1,3)--(1,3.5);
 \draw[thick] (1,-3)--(1,-3.5);

 \draw[thick, gray, dashed] (-2.25,0)--(4.5,0);

\draw [decorate,decoration={brace,amplitude=10pt}] (5,4) -- (5,0.2) node[midway,xshift=1cm]{$T_+$};
\draw [decorate,decoration={brace,amplitude=10pt}] (5,-0.2) -- (5,-4) node[midway,xshift=1cm]{$T_-$};

 \end{tikzpicture}
 \]
 \[
\begin{tikzpicture}[x=.35cm, y=.35cm,
    every edge/.style={
        draw,
      postaction={decorate,
                    decoration={markings}
                   }
        }
]

\node at (-3.5,0) {$\scalebox{1}{$T_+=$}$};
\node at (-1.25,-3.25) {\;};

\draw[thick] (2,2)--(1,3)--(-2,0);

\draw[thick] (0,0) -- (2,2)--(4,0);
 \draw[thick] (1,1) -- (2,0);

 \draw[thick] (1,3)--(1,3.5);

          \draw[fill=red] (0,0) circle (1pt); 
        \draw[fill=red] (2,0) circle (1pt); 
        \draw[fill=red] (4,0) circle (1pt); 
        \draw[fill=red] (-2,0) circle (1pt); 
        \draw[fill=green] (1,3) circle (1pt); 
        \draw[fill=black] (1,3.5) circle (1pt); 

 \end{tikzpicture}\qquad
\begin{tikzpicture}[x=.35cm, y=.35cm,
    every edge/.style={
        draw,
      postaction={decorate,
                    decoration={markings}
                   }
        }
]

\node at (-3.5,0) {$\scalebox{1}{$T_-=$}$};
\node at (-1.25,-3.25) {\;};

\draw[thick] (2,2)--(1,3)--(-2,0);

\draw[thick] (2,2)--(4,0);
 \draw[thick]  (2,2)--(0,0);
 \draw[thick] (3,1) -- (2,0);

 \draw[thick] (1,3)--(1,3.5);
  
          \draw[fill=red] (0,0) circle (1pt); 
        \draw[fill=red] (2,0) circle (1pt); 
        \draw[fill=red] (4,0) circle (1pt); 
        \draw[fill=red] (-2,0) circle (1pt); 
        \draw[fill=green] (1,3) circle (1pt); 
        \draw[fill=black] (1,3.5) circle (1pt); 

 \end{tikzpicture}
\]
\caption{
Top row: the two generators $x_0$ and $x_1$ of $F = F_2$ are shown. The dashed gray line is the $x$-axis; its intersection with the tree diagrams gives the leaves.\\
Bottom row: the top and bottom trees of $x_1$ are displayed. The leaves of $T_+$ and $T_-$ are the degree‑one vertices (circled in red).  
The root of $T_+$ and $T_-$ 
are circled in green. 
There are two other vertices of degree one, one in  
$T_+$ and one in $T_-$, they are circled in black.
}\label{genThompsonF}
\end{figure}
An alternative description of the elements of $F$ relevant to this paper is as follows: they can be visualized as pairs $(T_+,T_-)$ of planar binary rooted trees, each with the same number of leaves. We call such pairs binary tree diagrams, or just tree diagrams.
Usually one tree is positioned upside down on top of the other, with $T_+$ as the top tree and $T_-$ as the bottom tree (see Figure \ref{genThompsonF}, where $T_+$ and $T_-$ for the generator $x_1$ are indicated), and the leaves of 
$T_+$ are identified with those of $T_-$.
Recall that the degree of a vertex is the number of edges incident to it: the leaves have degree $1$ when we consider the top and bottom trees separately (as in the bottom row of Figure \ref{genThompsonF}) and degree $2$ when seen in a tree diagram (see the top row of the same figure).
In the representation
of the elements of $F$ as tree diagrams, there are only two vertices of degree $1$: 
  the highest one and the lowest one (see Figure \ref{genThompsonF}).

\begin{convention}\label{conventiondrawings}
We think of the tree diagrams being drawn on the $xy$-plane. 
The leaves sit on the $x$-axis, precisely  they
 are located at the points of $\mathbb{N}_0 :=\{0, 1, 2, \ldots \}$.
The only two  vertices of degree one are positioned respectively on the lines $y=1$ and $y=-1$.  
 \end{convention}

Two pairs of tree diagrams are regarded as equivalent if they differ only by pairs of opposing carets. In other words, two tree diagrams are equivalent if one can be transformed into the other by cancelling or inserting
 pairs of carets that are positioned opposite to each other. Figure \ref{cancellingcarets} shows an example of
 two equivalent tree diagrams differing by the cancellation/insertion of a pair of opposing carets.

\begin{figure}
\phantom{This text will be invisible} 
\[\begin{tikzpicture}[x=.5cm, y=.5cm,
    every edge/.style={
        draw,
      postaction={decorate,
                    decoration={markings}
                   }
        }
]

 \draw[thick] (0,0)--(1,1)--(2,0)--(1,-1)--(0,0); 
 \draw[thick] (1,1.5)--(1,1); 
 \draw[thick] (1,-1.5)--(1,-1); 
\node at (0,-1.2) {$\;$};
\node at (3.5,0) {$\scalebox{1}{$\leftrightarrow$}$};

\end{tikzpicture}
\begin{tikzpicture}[x=.5cm, y=.5cm,
    every edge/.style={
        draw,
      postaction={decorate,
                    decoration={markings}
                   }
        }
]

  \draw[thick] (1,1.5)--(1,-1.5); 
\node at (0,-1.2) {$\;$};
 
\end{tikzpicture}\qquad
\;\;
\begin{tikzpicture}[x=.35cm, y=.35cm,
    every edge/.style={
        draw,
      postaction={decorate,
                    decoration={markings}
                   }
        }
]

\node at (5.7,-1.5) {$\scalebox{1}{$\thicksim$}$};
\node at (-1.25,-3.25) {\;};

\begin{scope}[yshift=-.5cm, scale=1.5]  

\draw[thick] (0,0) -- (1.5,1.5) -- (3,0) -- (1.5,-1.5) -- (0,0); 
\draw[thick] (1,1) -- (2,0); 
\draw[thick] (1,0) -- (1.5,0.5); 
\draw[thick] (1,0) -- (2,-1); 
\draw[thick] (2,0) -- (1.5,-0.5);

\draw[thick] (1.5,2) -- (1.5,1.5); 
\draw[thick] (1.5,-1.5) -- (1.5,-2); 

 \end{scope}
 
\begin{scope}[xshift=2.5cm,yshift=-.5cm, scale=1]  

\draw[thick] (0,0) -- (2,2)--(4,0)--(2,-2)--(0,0);
 \draw[thick] (1,1) -- (2,0)--(3,-1);

 \draw[thick] (2,2)--(2,2.5);

 \draw[thick] (2,-2)--(2,-2.5);

 \end{scope}

 \end{tikzpicture}
\]
\caption{A pair  of cancelling carets and two equivalent tree diagrams.
}\label{cancellingcarets}
\end{figure}
The correspondence between homeomorphisms and tree diagrams is described in \cite[Section 2]{CFP}, which we recall briefly.
First, a standard dyadic interval is one of the form $[a/2^k, (a+1)/2^k]$, where $a, k\in\IZ$.
A standard dyadic partition of $[0,1]$ is a finite collection of standard dyadic intervals whose union is $[0,1]$, and 
where the intervals may overlap 
at most on a single point. These intervals and hence the partitions are in a one-to-one correspondence with the 
leaves of finite subtrees of the infinite binary planar tree $\CT$ (its root represents $[0,1]$, its children $[0,1/2]$ and $[1/2,1]$, etc).
It can be proved that for each $g\in F$ there exists a  standard dyadic partition of $[0,1]$ such that $g$ is differentiable on the interior of each interval of the partition, and the images of the intervals form   a standard dyadic partition of $[0,1]$ (see \cite[Lemma 2.2]{CFP}). 
This allows one to translate the description of the homeomorphisms in $F$ into the language of pairs of trees. This correspondence is not injective because the partition we started with can be refined, i.e., one can divide any of the intervals into two halves.

Moving forward, we consider a family of groups that generalizes Thompson’s group $F$: the Brown-Thompson groups. For any integer $p \geq 2$, the Brown-Thompson group $F_p$ is defined by the following presentation
\[
\langle x_0, x_1, \ldots \; | \; x_nx_k = x_kx_{n+p-1} \quad \forall \; k<n \rangle.
\]
The elements $x_0, x_1, \dots, x_{p-1}$ generate $F_p$, and the neutral element is denoted by $1$.

 \begin{figure}
\phantom{This text will be invisible} \[
\begin{tikzpicture}[x=.35cm, y=.35cm,
    every edge/.style={
        draw,
      postaction={decorate,
                    decoration={markings}
                   }
        }
]

\node at (-1.5,0) {$\scalebox{1}{$x_0=$}$};
\node at (-1.25,-3) {\;};

\draw[thick] (0,0) -- (2,2)--(4,0)--(2,-2)--(0,0);
\draw[thick] (1,1) -- (1,0)--(2,-2);
\draw[thick] (1,1) -- (2,0)--(3,-1);
\draw[thick] (2,2) -- (3,0)--(3,-1);

 \draw[thick] (2,2)--(2,2.5);
 \draw[thick] (2,-2)--(2,-2.5);

\end{tikzpicture}
\;\;
\begin{tikzpicture}[x=.35cm, y=.35cm,
    every edge/.style={
        draw,
      postaction={decorate,
                    decoration={markings}
                   }
        }
]

\node at (-1.5,0) {$\scalebox{1}{$x_1=$}$};
\node at (-1.25,-3.25) {\;};


\draw[thick] (0,0) -- (2,2)--(4,0)--(2,-2)--(0,0);
 \draw[thick] (1,0)--(2,-2);
\draw[thick] (3,0)--(2,1) -- (1,0); 
\draw[thick] (2,0)--(3,-1);
\draw[thick] (2,2) -- (2,0);
\draw[thick] (3,0)--(3,-1);


 \draw[thick] (2,2)--(2,2.5);
 \draw[thick] (2,-2)--(2,-2.5);
\end{tikzpicture}
\;\;
\begin{tikzpicture}[x=.35cm, y=.35cm,
    every edge/.style={
        draw,
      postaction={decorate,
                    decoration={markings}
                   }
        }
]

\node at (-3.5,0) {$\scalebox{1}{$x_2=$}$};
\node at (-1.25,-3.25) {\;};

\draw[thick] (2,2)--(1,3)--(-2,0)--(1,-3)--(2,-2);

\draw[thick] (0,0) -- (2,2)--(4,0)--(2,-2)--(0,0);
 \draw[thick] (1,1) -- (2,0)--(3,-1);

\draw[thick] (0,0) -- (2,2)--(4,0)--(2,-2)--(0,0);
\draw[thick] (1,1) -- (1,0)--(2,-2);
\draw[thick] (1,1) -- (2,0)--(3,-1);
\draw[thick] (2,2) -- (3,0)--(3,-1);

\draw[thick] (1,3) -- (-1,0)--(1,-3);

\node at (6,0) {$\ldots$};

 \draw[thick] (1,3)--(1,3.5);
 \draw[thick] (1,-3)--(1,-3.5);

\end{tikzpicture}
\]
\caption{The generators of the Brown-Thompson group $F_3$.}\label{genThompsonF3}
\end{figure}

As for $F$, the elements of $F_p$ can be described using pairs of $p$-ary trees. An illustration of the generators of $F_3$ is presented in Figure \ref{genThompsonF3}.  

\begin{convention}\label{conventiondrawings2}
As for $F$, 
we establish a convention for the visualization of trees in the plane in case of elements of $F_p$. 
The tree diagrams are thought as being drawn on the $xy$-plane, with 
the leaves being on the $x$-axis, 
more precisely  they
 are located at the points of $\mathbb{N}_0 :=\{0, 1, 2, \ldots \}$.
The only two vertices of degree one are positioned respectively on the lines $y=1$ and $y=-1$.  
\end{convention}

A natural embedding $\iota: F=F_2\to F_3$ is obtained by replacing each $3$-valent vertex in a tree diagram of $F_2$ with a $4$-valent vertex and connecting the middle edges in the only possible planar way. Denote temporarily by $\{y_i\}_{i\geq 0}$ the generators of $F$ and by $\{x_i\}_{i\geq 0}$ the generators of $F_3$. Under this embedding, we observe that $\iota(y_i)=x_{2i}$. For simplicity, we will often omit the symbol $\iota$ in our notation.

Let $(T_+,T_-)$ represent an element of $F_3$. We now provide a brief recapitulation, following \cite{Jo14, Jo18}, of how to construct $\Gamma(T_+,T_-)$, referred to as the \emph{planar graph of $(T_+,T_-)$}, \cite[Section 4.1]{Jo14}.
An  example of the whole construction is given in Figure \ref{Gammagraph}.

Each tree diagram $(T_+,T_-)$ should be thought of in the $xy$-plane, as explained in Convention \ref{conventiondrawings2}, with the two lines $y=1$ and $y=-1$ enclosing it. These two lines yield a strip containing the tree diagram. We color the strip alternately in black and white, with the leftmost region conventionally black. The vertices of the planar graph $\Gamma(T_+,T_-)$ are positioned along the $x$-axis at
\[
-1/2+2\mathbb{N}_0 :=\{-1/2, 1+1/2, 3+1/2, \ldots \},
\]
with one vertex corresponding to each black region. An edge is drawn between two black regions whenever they meet at a $4$-valent vertex. 
\begin{figure}
\[
\begin{tikzpicture}[x=.35cm, y=.35cm,
    every edge/.style={
        draw,
      postaction={decorate,
                    decoration={markings}
                   }
        }
]

 \node at (-1.25,-3.25) {\;};

\draw[thick] (2,2)--(1,3)--(-2,0)--(1,-3)--(2,-2);

\draw[thick] (0,0) -- (2,2)--(4,0)--(2,-2)--(0,0);
 \draw[thick] (1,1) -- (2,0)--(3,-1);

\draw[thick] (0,0) -- (2,2)--(4,0)--(2,-2)--(0,0);
\draw[thick] (1,1) -- (1,0)--(2,-2);
\draw[thick] (1,1) -- (2,0)--(3,-1);
\draw[thick] (2,2) -- (3,0)--(3,-1);

\draw[thick] (1,3) -- (-1,0)--(1,-3);

 \draw[thick] (1,3)--(1,3.5);
 \draw[thick] (1,-3)--(1,-3.5);

 \draw[thick] (-4,-3.5)--(6,-3.5);
 \draw[thick] (-4,3.5)--(6,3.5);

\end{tikzpicture}
\qquad
\begin{tikzpicture}[x=.35cm, y=.35cm,
    every edge/.style={
        draw,
      postaction={decorate,
                    decoration={markings}
                   }
        }
]

 \node at (-1.25,-3.25) {\;};

\draw[thick] (2,2)--(1,3)--(-2,0)--(1,-3)--(2,-2);

\draw[thick] (0,0) -- (2,2)--(4,0)--(2,-2)--(0,0);
 \draw[thick] (1,1) -- (2,0)--(3,-1);

\draw[thick] (0,0) -- (2,2)--(4,0)--(2,-2)--(0,0);
\draw[thick] (1,1) -- (1,0)--(2,-2);
\draw[thick] (1,1) -- (2,0)--(3,-1);
\draw[thick] (2,2) -- (3,0)--(3,-1);

\draw[thick] (1,3) -- (-1,0)--(1,-3);

 \draw[thick] (1,3)--(1,3.5);
 \draw[thick] (1,-3)--(1,-3.5);

 \draw[thick] (-4,-3.5)--(6,-3.5);
 \draw[thick] (-4,3.5)--(6,3.5);

\fill[gray!20] (-4,3.4)--(.9,3.4) -- (.9,3)--(-2.2,0) --  (0.9,-3) -- (0.9,-3.4) -- (-4,-3.4)-- cycle;

\fill[gray!20] (-.9,0) -- (1,2.9) -- (1.9,2) -- (-.1,0) -- (1.9,-2)-- (1, -2.9) -- cycle;

\fill[gray!20] (1.05,0)--(1.07,.85)--(2.9,-1) --(2,-1.9) -- cycle;

\fill[gray!20] (3.1,0)--(3.1,-.8)--(3.9,0)--(2.2, 1.72) -- cycle;

\end{tikzpicture}
\]
\[
\begin{tikzpicture}[x=.35cm, y=.35cm,
    every edge/.style={
        draw,
      postaction={decorate,
                    decoration={markings}
                   }
        }
]

 \node at (-1.25,-3.25) {\;};

\draw[thick] (2,2)--(1,3)--(-2,0)--(1,-3)--(2,-2);

\draw[thick] (0,0) -- (2,2)--(4,0)--(2,-2)--(0,0);
 \draw[thick] (1,1) -- (2,0)--(3,-1);

\draw[thick] (0,0) -- (2,2)--(4,0)--(2,-2)--(0,0);
\draw[thick] (1,1) -- (1,0)--(2,-2);
\draw[thick] (1,1) -- (2,0)--(3,-1);
\draw[thick] (2,2) -- (3,0)--(3,-1);

\draw[thick] (1,3) -- (-1,0)--(1,-3);

 \draw[thick] (1,3)--(1,3.5);
 \draw[thick] (1,-3)--(1,-3.5);

 \draw[thick] (-4,-3.5)--(6,-3.5);
 \draw[thick] (-4,3.5)--(6,3.5);

\fill[gray!20] (-4,3.4)--(.9,3.4) -- (.9,3)--(-2.2,0) --  (0.9,-3) -- (0.9,-3.4) -- (-4,-3.4)-- cycle;

\fill[gray!20] (-.9,0) -- (1,2.9) -- (1.9,2) -- (-.1,0) -- (1.9,-2)-- (1, -2.9) -- cycle;

\fill[gray!20] (1.05,0)--(1.07,.85)--(2.9,-1) --(2,-1.9) -- cycle;

\fill[gray!20] (3.1,0)--(3.1,-.8)--(3.9,0)--(2.2, 1.72) -- cycle;

\draw[fill=black] (-2.5,0) circle (1pt); 
\draw[fill=black] (-.5,0) circle (1pt); 
\draw[fill=black] (1.5,0) circle (1pt); 
\draw[fill=black] (3.5,0) circle (1pt); 
\end{tikzpicture}
\qquad
\qquad
\begin{tikzpicture}
\begin{scope}[yshift=3cm, scale=.75]
    \draw (0,0) to[out=90, in=90] (1,0);  
    \draw (1,0) to[out=90, in=90] (2,0);  
    \draw (1,0) to[out=90, in=90] (3,0);  
    
    \draw (0,0) to[out=-90, in=-90] (1,0);  
    \draw (1,0) to[out=-90, in=-90] (2,0);  
    \draw (2,0) to[out=-90, in=-90] (3,0);

         \draw[fill=black] (0,0) circle (2pt); 
        \draw[fill=black] (1,0) circle (2pt); 
        \draw[fill=black] (2,0) circle (2pt); 
        \draw[fill=black] (3,0) circle (2pt); 
\end{scope}        
\useasboundingbox (1,1.75) rectangle (5,5); 
\end{tikzpicture}
\]
\caption{Top row: the strip for $x_1^2 = (T_+,T_-) \in F_3$ and the corresponding shading. Bottom row: the vertices for the corresponding planar graph sitting in the shaded regions of the strip, and the planar graph $\Gamma(T_+,T_-)$.}
\label{Gammagraph}
\end{figure}

For an element $(T_+,T_-)$ in $F_2$, we define the \emph{planar graph of $(T_+,T_-)$} to be that of $\iota(T_+,T_-)\in F_3$.

We anticipate that the planar graph $\Gamma(T_+,T_-)$ is closely related to the Tait graph of the link $\mathcal{L}(T_+,T_-)$, except that we do not specify the signs of its edges. 
As for the graphs corresponding to elements of $F$, the signs of the edges are actually $+$ if they are in the upper-half plane, $-$ if they are in the lower-half plane. 
We will come back to this topic in Section \ref{sec_knots}.
See also \cite[Section 5.3.1]{Jo14}. 

Jones introduced the oriented subgroups in two works: the first in \cite{Jo14}, and the second in \cite{Jo18}. We are now in a position to define them:
\begin{align*}
\vec{F}&:=\{(T_+,T_-)\in F \; | \; \Gamma(T_+,T_-) \textrm{ is $2$-colorable} \} \\  
&:=\{(T_+,T_-)\in F \; | \; {\rm Chr}_{\Gamma(T_+,T_-)}(2)=2\}, \\  
\vec{F}_3&:=\{(T_+,T_-)\in F_3\; | \; \Gamma(T_+,T_-) \textrm{ is $2$-colorable} \}\\
&:=\{(T_+,T_-)\in F_3 \; | \; {\rm Chr}_{\Gamma(T_+,T_-)}(2)=2\}, 
\end{align*}
where ${\rm Chr}_{\Gamma}(t)$ denotes the chromatic polynomial. A graph is said to be \emph{$2$-colorable} if its vertices can be labeled with two colors such that adjacent vertices have different colors. We designate these colors as $+$ and $-$.

Importantly, 
since $\Gamma(T_+,T_-)$ is connected, if it is $2$-colorable, there exist precisely two valid colorings. By convention, we choose the one in which the leftmost vertex is assigned the color $+$. Notably, for $\Gamma(T_+,T_-)$, being $2$-colorable is equivalent to being bipartite.

We now briefly recall Jones' construction of knots from $\vec{F}$.
First we turn all the $3$-valent vertices into $4$-valent by means of the map $\iota$.
Then we draw an edge between the two roots of the trees. 
The third step consists of replacing all $4$-valent vertices 
by the forks  in Figure \ref{rulesknot}.
We exemplify this procedure with $x_0x_1$ in Figure \ref{knotproducer}. 
\begin{figure}
\[
\begin{tikzpicture}[x=.3cm, y=.3cm,
    every edge/.style={
        draw,
      postaction={decorate,
                    decoration={markings}
                   }
        }
]

\draw[thick] (1,1)--(1,2);
\draw[thick] (0,0) --(1,1)--(2,0);

\node at (0,-1.2) {$\;$};

\end{tikzpicture}\quad
\begin{tikzpicture}[x=.3cm, y=.3cm,
    every edge/.style={
        draw,
      postaction={decorate,
                    decoration={markings}
                   }
        }
]

\draw[thick] (1,0)--(1,2);
\draw[thick] (0,0) --(1,1)--(2,0);

\node at (0,-1.2) {$\;$};
\node at (-2,1) {$\scalebox{1}{$\mapsto\; $}$};

\end{tikzpicture}\qquad
\begin{tikzpicture}[x=.3cm, y=.3cm,
    every edge/.style={
        draw,
      postaction={decorate,
                    decoration={markings}
                   }
        }
]

\draw[thick] (1,0)--(1,2);
\draw[thick] (0,0) --(1,1)--(2,0);

\node at (0,-1.2) {$\;$};

\end{tikzpicture}\quad
\begin{tikzpicture}[x=.3cm, y=.3cm,
    every edge/.style={
        draw,
      postaction={decorate,
                    decoration={markings}
                   }
        }
]

\draw[thick] (1,0)--(1,.35);
\draw[thick] (1,.75)--(1,2);
\draw[thick] (0,0) to[out=90,in=90] (2,0);

\node at (-2,1) {$\scalebox{1}{$\mapsto\; $}$};

\node at (0,-1.2) {$\;$};

\end{tikzpicture}
\qquad \begin{tikzpicture}[x=.3cm, y=.3cm,
    every edge/.style={
        draw,
      postaction={decorate,
                    decoration={markings}
                   }
        }
]

\draw[thick] (1,0)--(1,2);
\draw[thick] (0,2) --(1,1)--(2,2);

\node at (0,-1.2) {$\;$};

\end{tikzpicture}\quad
\begin{tikzpicture}[x=.3cm, y=.3cm,
    every edge/.style={
        draw,
      postaction={decorate,
                    decoration={markings}
                   }
        }
]

\draw[thick] (1,2)--(1,1.65);
\draw[thick] (1,1.25)--(1,0);
\draw[thick] (0,2) to[out=-90,in=-90] (2,2);

\node at (-2,1) {$\scalebox{1}{$\mapsto\; $}$};

\node at (0,-1.2) {$\;$};

\end{tikzpicture}
\]
\caption{The rules needed for obtaining $\CL(g)$. }\label{rulesknot}
\end{figure}

\begin{figure}
\[
\begin{tikzpicture}[x=.35cm, y=.35cm,
    every edge/.style={
        draw,
      postaction={decorate,
                    decoration={markings}
                   }
        }
]
 
\draw[thick] (0,0) -- (3,3)--(6,0)--(3,-3)--(0,0);
 \draw[thick] (3,1) -- (2,0)--(4,-2);

 \draw[thick] (2,2)--(5,-1);

 \draw[thick] (3,-3)--(3,-3.5);
 \draw[thick] (3,3)--(3,3.5);

\node at (7,0) {$\scalebox{1}{$\mapsto$}$};

\node at (0,-1.2) {\phantom{$\frac{T_+}{T_-}=$}};
\node at (0,-3) {\phantom{$\frac{T_+}{T_-}=$}};

\end{tikzpicture}
\begin{tikzpicture}[x=.35cm, y=.35cm,
    every edge/.style={
        draw,
      postaction={decorate,
                    decoration={markings}
                   }
        }
]

 \draw[thick] (3,3.5) -- (3,3);  
 \draw[thick] (3,-3.5) -- (3,-3);  

\draw[thick] (0,0) -- (3,3)--(6,0)--(3,-3)--(0,0);
 \draw[thick] (3,1) -- (2,0)--(4,-2);

 \draw[thick] (2,2)--(5,-1);

\node at (7,0) {$\scalebox{1}{$\mapsto$}$};

 \draw[thick] (2,2)--(1,0)--(3,-3);
 \draw[thick] (3,1)--(3,0)--(4,-2);
 \draw[thick] (5,-1)--(5,0)--(3,3);

\node at (0,-1.2) {\phantom{$\frac{T_+}{T_-}=$}};
\node at (0,-3) {\phantom{$\frac{T_+}{T_-}=$}};

\end{tikzpicture}
\begin{tikzpicture}[x=.35cm, y=.35cm,
    every edge/.style={
        draw,
      postaction={decorate,
                    decoration={markings}
                   }
        }
]

 \draw[thick] (-1,3) to[out=90,in=90] (3,3);  
 \draw[thick] (-1,-3) to[out=-90,in=-90] (3,-3);  
\draw[thick] (-1,-3) --(-1,3);

\draw[thick] (0,0) -- (3,3)--(6,0)--(3,-3)--(0,0);
 \draw[thick] (3,1) -- (2,0)--(4,-2);

 \draw[thick] (2,2)--(5,-1);

\node at (7,0) {$\scalebox{1}{$\mapsto$}$};

 \draw[thick] (2,2)--(1,0)--(3,-3);
 \draw[thick] (3,1)--(3,0)--(4,-2);
 \draw[thick] (5,-1)--(5,0)--(3,3);

\node at (0,-1.2) {\phantom{$\frac{T_+}{T_-}=$}};
\node at (0,-3) {\phantom{$\frac{T_+}{T_-}=$}};

\end{tikzpicture}
\begin{tikzpicture}[x=.35cm, y=.35cm,
    every edge/.style={
        draw,
      postaction={decorate,
                    decoration={markings}
                   }
        }
]

\draw[thick] (2,0) to[out=90,in=90] (4,0);   
\draw[thick] (0,0) to[out=90,in=90] (3,.75);   
\draw[thick] (1,1.25) to[out=90,in=90] (6,0);   
\draw[thick] (3,2.4) to[out=90,in=90] (-1,0);   
\draw[thick] (1,.9)--(1,0);
\draw[thick] (3,.4)--(3,0);
\draw[thick] (3,2.1) to[out=-90,in=90] (5,0);

\draw[thick] (4,0) to[out=-90,in=-90] (6,0);   
\draw[thick] (2,0) to[out=-90,in=-90] (5,-.8);   
\draw[thick] (0,0) to[out=-90,in=-90] (3,-1.3);   
\draw[thick] (-1,0) to[out=-90,in=-90] (2,-2);   
\draw[thick] (5,0)--(5,-.4);
\draw[thick] (3,0)--(3,-.9);
\draw[thick] (1,0) to[out=-90,in=90] (2,-1.6);
 
\node at (0,-4) {$\scalebox{1}{$\,$}$};

\node at (0,-1.2) {$\;$};

\end{tikzpicture}
\qquad
\]
\caption{From left to right: the element $x_0x_1$; its image in $F_3$ under the embedding map $\iota$; its closure; and the corresponding link.}
\label{knotproducer}
\end{figure}
The Tait graph of the resulting knot is exactly the graph $\Gamma(g)$ and, 
by using the coloring fixed by convention for $\vec{F}$, we obtain an oriented knot/link: if a vertex  with color $+$ sits inside a region, then its boundary is oriented counterclockwise, if the color is $-$ the orientation is clockwise. Given a tree diagram $(T_+,T_-)\in \vec{F}$, we denote by $\vec{\CL}(T_+,T_-)$ the corresponding knot/link.

We now pass to the  $3$-colorable subgroup $\CF$. 
As any tree diagram of $F$ is drawn on the $xy$-plane according to Convention \ref{conventiondrawings}, it sits on the strip bounded by the lines $y=1$ and $y=-1$. If we consider the complement of the tree diagram in this strip, we have a (finite) number 
of connected components, which we call the partition of the strip. Figure \ref{3colorabilityexample} shows an element of $F$ and its strip.

 \begin{figure}
\phantom{This text will be invisible} \[
\begin{tikzpicture}[x=.35cm, y=.35cm,
    every edge/.style={
        draw,
      postaction={decorate,
                    decoration={markings}
                   }
        }
]

 \node at (-1.25,-3) {\;};

\draw[thick] (0,0)--(2,2)--(4,0)--(2,-2)--(0,0);

\draw[thick] (1,0)--(1.5, .5);
\draw[thick] (2,0)--(1, 1);
\draw[thick] (3,0)--(1.5, 1.5);

\draw[thick] (1,0)--(2.5, -1.5);
\draw[thick] (2,0)--(3, -1);
\draw[thick] (3,0)--(2.5, -.5);

\draw[thick] (2,2)--(2,2.5);
\draw[thick] (2,-2)--(2,-2.5);

\end{tikzpicture}
\;\;
\qquad
\begin{tikzpicture}[x=.35cm, y=.35cm,
    every edge/.style={
        draw,
      postaction={decorate,
                    decoration={markings}
                   }
        }
]

 \node at (-1.25,-3) {\;};

\draw[thick] (0,0)--(2,2)--(4,0)--(2,-2)--(0,0);

\draw[thick] (1,0)--(1.5, .5);
\draw[thick] (2,0)--(1, 1);
\draw[thick] (3,0)--(1.5, 1.5);

\draw[thick] (1,0)--(2.5, -1.5);
\draw[thick] (2,0)--(3, -1);
\draw[thick] (3,0)--(2.5, -.5);

\draw[thick] (2,2)--(2,2.5);
\draw[thick] (2,-2)--(2,-2.5);

\draw[thick] (-2,2.5)--(6,2.5);
\draw[thick] (-2,-2.5)--(6,-2.5);

\end{tikzpicture}
\;\;
\qquad
\begin{tikzpicture}[x=.35cm, y=.35cm,
    every edge/.style={
        draw,
      postaction={decorate,
                    decoration={markings}
                   }
        }
]

 \node at (-1.25,-3) {\;};

\draw[thick] (0,0)--(2,2)--(4,0)--(2,-2)--(0,0);

\draw[thick] (1,0)--(1.5, .5);
\draw[thick] (2,0)--(1, 1);
\draw[thick] (3,0)--(1.5, 1.5);

\draw[thick] (1,0)--(2.5, -1.5);
\draw[thick] (2,0)--(3, -1);
\draw[thick] (3,0)--(2.5, -.5);

\draw[thick] (2,2)--(2,2.5);
\draw[thick] (2,-2)--(2,-2.5);

\draw[thick] (-2,2.5)--(6,2.5);
\draw[thick] (-2,-2.5)--(6,-2.5);

\node at (-1,0) {$\scalebox{.75}{$0$}$};
\node at (.75,0) {$\scalebox{.75}{$2$}$};
\node at (1.5,0) {$\scalebox{.75}{$0$}$};
\node at (2.5,0) {$\scalebox{.75}{$1$}$};
\node at (3.5,0) {$\scalebox{.75}{$2$}$};
\node at (5,0) {$\scalebox{.75}{$1$}$};
 
\end{tikzpicture}
\]
\caption{
    From left to right: an element of $F$ (which also lies in $\CF$); the corresponding strip on the $xy$-plane; and the unique admissible coloring according to rule \eqref{convention-colors}.
}\label{3colorabilityexample}
\end{figure}
This strip may or may not be $3$-colorable, meaning it may or may not be possible to assign colors from $\IZ_3 = \{0, 1, 2\}$ to its regions such that any two regions sharing an edge are assigned different colors. 
By convention, we fix specific colors for the regions near the roots. 
  \begin{eqnarray}\label{convention-colors}
\begin{tikzpicture}[x=.5cm, y=.5cm,
    every edge/.style={
        draw,
      postaction={decorate,
                    decoration={markings}
                   }
        }
] 

 \draw[thick] (0,0) -- (1,1)--(2,0);
 \draw[thick] (1,1) -- (1,2);



\node at (0,1) {$0$};
\node at (2,1) {$1$};
\node at (1,0.2) {$2$};

 \draw[thick] (-1,2) -- (3,2);


\end{tikzpicture}
\qquad\qquad 
\begin{tikzpicture}[x=.5cm, y=.5cm,
    every edge/.style={
        draw,
      postaction={decorate,
                    decoration={markings}
                   }
        }
] 

 \draw[thick] (-1,0) -- (3,0);

 \draw[thick] (0,2) -- (1,1)--(2,2);
 \draw[thick] (1,1) -- (1,0);



\node at (0,1) {$0$};
\node at (2,1) {$1$};
\node at (1,1.8) {$2$};

\end{tikzpicture}
\end{eqnarray}
Once this convention is applied, if the strip is $3$-colorable, the coloring is uniquely determined. 
The $3$-colorable subgroup $\CF$ consists of all elements of $F$ for which the corresponding strip admits such a $3$-coloring.
This subgroup is contained in the rectangular subgroup $K_{(2,2)}$, which belongs to the  
family of groups 
\begin{align} \label{rectsub22}
K_{(a,b)}:=\{g\in F \, | \, \log_2 g'(0)\in a\IZ, \log_2g'(1)\in b\IZ\}
\end{align}
$a, b\in\IN =\{1,2,3,\ldots\}$. This family was introduced in \cite{BW}.
All these groups have   finite index in $F$ and 
are isomorphic to it.
 Note that $\log_2 g'(0)$ can be read from the tree diagram as the difference between the path from the root of the top tree to the left-most leaf and  the path from the root of the bottom tree to the left-most leaf. Similarly,  $\log_2 g'(1)$ can be read from the tree diagram as the difference between the path from the root of the top tree to the right-most leaf
and  the path from the root of the bottom tree to the right-most leaf. Interestingly, all knots arising from $\CF$ are $3$-colorable, \cite{Kodama3}.

\section{The subgroups $M_0$, $M_1$, and $M_2$} \label{section_main_res}
Denote by $\sigma: F\to F$ the order-$2$ automorphism defined on the level of tree diagrams by a reflection about any vertical line.
In \cite[Section 4]{TV2}, three maximal subgroups of infinite index in $K_{(2,2)}$ were identified
\begin{align*}
& M_0:=\langle \CF, x_0^2\rangle, \qquad
 M_1:=\langle \CF, x_1^2\rangle, \qquad
 M_2:=\langle \CF, \sigma(x_1^2)\rangle.
\end{align*}
Recall that since $K_{(2,2)}$ is isomorphic to $F$, these subgroups yield maximal subgroups of infinite index in $F$.
For any $t\in [0,1]$, denote by $.a_1 \ldots a_n$ its dyadic expansion.
Now consider the subsets 
\[
S_i = \{ t \in (0,1) \cap \IZ[1/2] \mid \omega(t) \equiv_3 i \},
\]
where \( i \in \IZ_3 \), \( \equiv_3 \) denotes equivalence modulo \( 3 \), and the weight function \( \omega \) is defined by 
\[
\omega(.a_1 \ldots a_n) := \sum_{j=1}^n (-1)^j a_j \pmod{3}.
\]
These subsets \( S_i \) form a partition of the dyadic rationals in the interval \( (0,1) \). The $3$-colorable subgroup $\CF$ stabilizes these subsets\footnote{For this discussion we refer the reader to the arXiv version of \cite{TV3}, arXiv:2212.12269v4.}, more precisely $\CF=\cap_{i\in\IZ_3} {\rm Stab}(S_i)$ \cite[Lemma 2.7]{TV3}. 
The following theorem shows that the subgroups $M_0$, $M_1$, $M_2$ can be described as stabilizers of these subsets. 
\setcounter{theorem}{0}
\begin{theorem}\label{stabtheo}
The following equalities hold
\[
M_0 = \stab(S_2), \quad 
M_1 = \stab(S_0) \cap K_{(2,2)}, \quad
M_2 = \stab(S_1) \cap K_{(2,2)}.
\]
\end{theorem}
Before undertaking the proof, we make a couple of general observations that we will use repeatedly in this proof. 
\begin{remark}\label{stabilizing}
In general, we may color all the regions of the strip contained in the upper-half plane and in the lower-half plane separately, always following Convention \ref{convention-colors}.
If the colors of the two halves of the strip match, the element is in $\CF$. Otherwise,
fix $i\in\{0,1,2\}$.
If the regions with color $i$ match, then the element is in $\stab(S_i)$.
\end{remark}
\begin{remark}\label{color_at_extremes}
The color of 
the region to the left of the rightmost leaf in the top/bottom tree is determined separately by
each tree:
the color is  $2$ 
if the length of the path from the root to the leaf is odd, and $0$ otherwise. 

Similarly, the color for the
region to the right of the leftmost leaf is determined separately by the top tree and the bottom tree:
the color is $2$
if the path from the root to the first leaf has odd length, and $1$ otherwise.
\end{remark}
\begin{proof}[Proof of Theorem \ref{stabtheo}]
We start by proving the inclusion $M_0 \subseteq \stab(S_2)$.
As the $3$-colorable subgroup $\CF$  is contained in $\stab(S_2)$, we only need to prove that 
  $x_0^2$ belongs to $\stab(S_2)$.
Let us draw  the  strip corresponding to this element.
\begin{equation}\label{primoese}
\begin{tikzpicture}[x=.35cm, y=.35cm,
    every edge/.style={
        draw,
      postaction={decorate,
                    decoration={markings}
                   }
        }
]
 \node at (-2,0) {$\scalebox{1}{$x_0^2=$}$}; 
 \node at (6,0.5) {\;};

\draw[thick] (0,0) -- (3,3)--(6,0)--(3,-3)--(0,0);
 \draw[thick] (1,1) -- (4,-2);
 \draw[thick] (2,2) -- (5,-1);

 \draw[thick] (3,3)--(3,3.5);

 \draw[thick] (3,-3)--(3,-3.5);

\end{tikzpicture}\qquad
\begin{tikzpicture}[x=.35cm, y=.35cm,
    every edge/.style={
        draw,
      postaction={decorate,
                    decoration={markings}
                   }
        }
]
\node at (-.5,0.5) {$\scalebox{.75}{$0$}$};
\node at (1,0.5) {$\scalebox{.75}{$2$}$};
\node at (3,0.5) {$\scalebox{.75}{$1$}$};
\node at (5,0.5) {$\scalebox{.75}{$2$}$};
\node at (7,0.5) {$\scalebox{.75}{$1$}$};

\node at (-.5,-0.5) {$\scalebox{.75}{$0$}$};
\node at (1,-0.5) {$\scalebox{.75}{$2$}$};
\node at (3,-0.5) {$\scalebox{.75}{$0$}$};
\node at (5,-0.5) {$\scalebox{.75}{$2$}$};
\node at (7,-0.5) {$\scalebox{.75}{$1$}$};
\draw[dashed] (-1,0)--(7,0);
\node at (6,0.5) {\;};

\draw[thick] (0,0) -- (3,3)--(6,0)--(3,-3)--(0,0);
 \draw[thick] (1,1) -- (4,-2);
 \draw[thick] (2,2) -- (5,-1);

 \draw[thick] (3,3)--(3,3.5);

 \draw[thick] (3,-3)--(3,-3.5);

 \draw[thick] (7,-3.5)--(-1,-3.5);
 \draw[thick] (7,3.5)--(-1,3.5);

\end{tikzpicture}
\end{equation}
The regions in the top and bottom trees labeled with color $2$ match.  Therefore, we have $M_0 \subseteq \stab(S_2)$. \\
Additionally, we note that there exists a region of color $1$ in the upper half-plane  
connected to a region of color $0$ in the lower half-plane, so this element does not preserve $S_1$.

It was proven in \cite[Theorem 4.6]{TV2} that $M_0$ is maximal in $K_{(2,2)}$.
This means that in order to prove that
$M_0 = \stab(S_2)$,
 it suffices to show that 
$\stab(S_2)$ is a proper subgroup of $K_{(2,2)}$.  
Let $g$ be an element of $\stab(S_2)$. 
%
Consider the color of the region to the left of the rightmost leaf in the top/bottom tree. 
This is determined separately by each tree according to
 Remark \ref{color_at_extremes} and can be either $0$ or $2$ depending on the parity of the length of the path from the root to the leaf.   
Since $g$ belongs to $\stab(S_2)$, these
 must be the same, which implies that $\log_2 g'(1) \in 2\mathbb{Z}$.  \\
Similarly, by  Remark \ref{color_at_extremes}
the color for the
region to the right of the leftmost leaf can be either $1$ or $2$.
As the element $g \in \stab(S_2)$, the paths from the first leaf to the roots of the top and bottom tree 
must   have the same parity, leading to the conclusion that $\log_2 g'(0) \in 2\mathbb{Z}$.  \\
Finally, it remains to show that $\stab(S_2)$ is a proper subgroup of $K_{(2,2)}$. 
The element $h_1$ 
in \eqref{eseh1}, 
 belongs to $K_{(2,2)}\setminus \stab(S_2)$. Indeed it is $K_{(2,2)}$
 because the paths from the leftmost leaf to the root have the same parity in both the top and bottom trees, and likewise for the rightmost leaf. It does not belong to $\stab(S_2)$ because there is a region colored $2$ on the top tree matching one colored $1$ on the bottom tree.
\begin{equation}\label{eseh1}
\begin{tikzpicture}[x=.6cm, y=.6cm,
    every edge/.style={
        draw,
      postaction={decorate,
                    decoration={markings}
                   }
        }
] 
\node at (6,0.5) {\;};

 \draw[thick] (0,0)--(2.5, 2.5)--(5,0)--(2.5, -2.5)--(0,0) -- (1.5, 1.5)--(3,0)--(2,-1)--(1,0)--(2,1);
 \draw[thick] (1.5, -.5) -- (2.5, .5);
 \draw[thick] (2,-1) -- (3, -2);
 \draw[thick] (2,2) -- (4.5, -.5);
 
 \draw[thick] (2.5,2.5)--(2.5,3);
 \draw[thick] (2.5,-2.5)--(2.5,-3);

\node at (-1,0) {$\scalebox{1}{$h_1=$}$};

\end{tikzpicture}\qquad
\begin{tikzpicture}[x=.6cm, y=.6cm,
    every edge/.style={
        draw,
      postaction={decorate,
                    decoration={markings}
                   }
        }
]
\node at (-.5,0.5) {$\scalebox{.6}{$0$}$};
\node at (1,0.5) {$\scalebox{.6}{$2$}$};
\node at (3,0.5) {$\scalebox{.6}{$1$}$};
\node at (5,0.5) {$\scalebox{.6}{$1$}$};
\node at (4,0.5) {$\scalebox{.6}{$2$}$};
\node at (2,0.5) {$\scalebox{.6}{$0$}$};
\node at (2.5,0.2) {$\scalebox{.6}{$2$}$};
 
\node at (-.5,-0.5) {$\scalebox{.6}{$0$}$};
\node at (1,-0.5) {$\scalebox{.6}{$2$}$};
\node at (3,-0.5) {$\scalebox{.6}{$0$}$};
\node at (5,-0.5) {$\scalebox{.6}{$1$}$};
\node at (4.5,-0.25) {$\scalebox{.6}{$2$}$};
\node at (2.5,-0.2) {$\scalebox{.6}{$1$}$};
\node at (1.5,-0.25) {$\scalebox{.6}{$0$}$};
 \draw[dashed] (-1,0)--(6,0);
\node at (6,0.5) {\;};

 \draw[thick] (0,0)--(2.5, 2.5)--(5,0)--(2.5, -2.5)--(0,0) -- (1.5, 1.5)--(3,0)--(2,-1)--(1,0)--(2,1);
 \draw[thick] (1.5, -.5) -- (2.5, .5);
 \draw[thick] (2,-1) -- (3, -2);
 \draw[thick] (2,2) -- (4.5, -.5);
 
 \draw[thick] (2.5,2.5)--(2.5,3);
 \draw[thick] (2.5,-2.5)--(2.5,-3);

 \draw[thick] (6,-3)--(-1,-3);
 \draw[thick] (6,3)--(-1,3);

\end{tikzpicture}
\end{equation}

We now turn to $M_1$ and we adopt the same strategy as for $M_0$. 
Here, we show that $x_1^2$ belongs to $\stab(S_0)$.  
\[
\begin{tikzpicture}[x=.35cm, y=.35cm,
    every edge/.style={
        draw,
      postaction={decorate,
                    decoration={markings}
                   }
        }
]

\node at (-4,0)  {$\scalebox{1}{$x_1^2=$}$}; 
 \node at (6,0.5) {\;};

\draw[thick] (0,0) -- (3,3)--(6,0)--(3,-3)--(0,0);
 \draw[thick] (1,1) -- (4,-2);
 \draw[thick] (2,2) -- (5,-1);

 \draw[thick] (3,3)--(2,4)--(-2,0)--(2, -4)--(3,-3);

 \draw[thick] (2,4)--(2,4.5);

 \draw[thick] (2,-4)--(2,-4.5);

\end{tikzpicture}\qquad
\begin{tikzpicture}[x=.35cm, y=.35cm,
    every edge/.style={
        draw,
      postaction={decorate,
                    decoration={markings}
                   }
        }
]
\node at (-2.5,0.5) {$\scalebox{.75}{$0$}$};
\node at (-.5,0.5) {$\scalebox{.75}{$2$}$};
\node at (1,0.5) {$\scalebox{.75}{$0$}$};
\node at (3,0.5) {$\scalebox{.75}{$1$}$};
\node at (5,0.5) {$\scalebox{.75}{$0$}$};
\node at (7,0.5) {$\scalebox{.75}{$1$}$};

\node at (-2.5,-0.5) {$\scalebox{.75}{$0$}$};
\node at (-.5,-0.5) {$\scalebox{.75}{$2$}$};
\node at (1,-0.5) {$\scalebox{.75}{$0$}$};
\node at (3,-0.5) {$\scalebox{.75}{$2$}$};
\node at (5,-0.5) {$\scalebox{.75}{$0$}$};
\node at (7,-0.5) {$\scalebox{.75}{$1$}$};
\draw[dashed] (-3,0)--(7,0);
\node at (6,0.5) {\;};

\draw[thick] (0,0) -- (3,3)--(6,0)--(3,-3)--(0,0);
 \draw[thick] (1,1) -- (4,-2);
 \draw[thick] (2,2) -- (5,-1);

 \draw[thick] (3,3)--(2,4)--(-2,0)--(2, -4)--(3,-3);

 \draw[thick] (2,4)--(2,4.5);

 \draw[thick] (2,-4)--(2,-4.5);

 \draw[thick] (7,-4.5)--(-3,-4.5);
 \draw[thick] (7,4.5)--(-3,4.5);

\end{tikzpicture}
\]
The subgroup $M_1$ is generated by $x_1^2$ and $\CF$. The latter is contained in $\stab(S_0)$ and, thus, $M_1$ is a subgroup of $\stab(S_0)$. It was proven in \cite[Theorem 4.6]{TV2} that $M_1$ is maximal in $K_{(2,2)}$. As   observed
in Remark \ref{color_at_extremes}, the region to the left of the rightmost leaf in the top/bottom tree is colored $2$ if the path from the root to the leaf has odd length, and $0$ otherwise. Since $g$ belongs to $\stab(S_0)$, these paths must have the same parity, which implies that $\log_2 g'(1) \in 2\mathbb{Z}$. Consequently, we obtain $\stab(S_0) \subseteq K_{(1,2)}$.
Unlike the first case, $\stab(S_0)$ is not contained in $K_{(2,2)}$. In fact, the subgroup $\stab(S_0)$ is strictly larger than $K_{(2,2)}$ as the element below belongs to the former but not the latter.
\[
x_0 x_1 \in \stab(S_0) \cap (K_{(1,2)} \setminus K_{(2,2)}).
\]
\[
\begin{tikzpicture}[x=.35cm, y=.35cm,
    every edge/.style={
        draw,
      postaction={decorate,
                    decoration={markings}
                   }
        }
] 
\node at (-4.5,0) {$\scalebox{1}{$x_0x_1=$}$}; 
\node at (6,0.5) {\;};

\draw[thick] (2,2)--(1,3)--(-2,0)--(1,-3)--(2,-2);

\draw[thick] (0,0) -- (1,1); 
\draw[thick] (0,2) -- (1,1); 
\draw[thick] (2,2)--(4,0)--(2,-2)--(0,0);
 \draw[thick] (1,1) -- (2,0)--(3,-1);

 \draw[thick] (1,3)--(1,3.5);
 \draw[thick] (1,-3)--(1,-3.5);

\end{tikzpicture}
\qquad
\begin{tikzpicture}[x=.35cm, y=.35cm,
    every edge/.style={
        draw,
      postaction={decorate,
                    decoration={markings}
                   }
        }
]
\node at (-2.5,0.5) {$\scalebox{.75}{$0$}$};
\node at (-.5,0.5) {$\scalebox{.75}{$1$}$};
\node at (1,0.5) {$\scalebox{.75}{$0$}$};
\node at (3,0.5) {$\scalebox{.75}{$2$}$};
\node at (5,0.5) {$\scalebox{.75}{$1$}$};
 
\node at (-2.5,-0.5) {$\scalebox{.75}{$0$}$};
\node at (-.5,-0.5) {$\scalebox{.75}{$2$}$};
\node at (1,-0.5) {$\scalebox{.75}{$0$}$};
\node at (3,-0.5) {$\scalebox{.75}{$2$}$};
\node at (5,-0.5) {$\scalebox{.75}{$1$}$};
 \draw[dashed] (-3,0)--(5,0);
\node at (6,0.5) {\;};

\draw[thick] (2,2)--(1,3)--(-2,0)--(1,-3)--(2,-2);

\draw[thick] (0,0) -- (1,1); 
\draw[thick] (0,2) -- (1,1); 
\draw[thick] (2,2)--(4,0)--(2,-2)--(0,0);
 \draw[thick] (1,1) -- (2,0)--(3,-1);

 \draw[thick] (1,3)--(1,3.5);
 \draw[thick] (1,-3)--(1,-3.5);

 \draw[thick] (5,-3.5)--(-3,-3.5);
 \draw[thick] (5,3.5)--(-3,3.5);

\end{tikzpicture}
\]
We are thus led to consider the subgroup $\stab(S_0) \cap K_{(2,2)}$. This subgroup contains $M_1$. 
We need to show that $\stab(S_0) \cap K_{(2,2)}$ is a proper subgroup of $K_{(2,2)}$, because then the claim follows from the maximality of $M_1$ in $K_{(2,2)}$. 
To this end, we may resort to the element $k_1$, 
displayed in \eqref{-eseh1}.
This element is not in $\stab(S_0)$ by Remark \ref{stabilizing}.
\begin{equation}\label{-eseh1}
\begin{tikzpicture}[x=.6cm, y=.6cm,
    every edge/.style={
        draw,
      postaction={decorate,
                    decoration={markings}
                   }
        }
] 
\node at (6,0.5) {\;};

 \draw[thick] (0,0)--(2.5, 2.5)--(5,0)--(2.5, -2.5)--(0,0) -- (1.5, 1.5)--(3,0)--(2,-1)--(1,0)--(2,1);
 \draw[thick] (1.5, .5) -- (2.5, -.5);
 \draw[thick] (2,-1) -- (3, -2);
 \draw[thick] (2,2) -- (4.5, -.5);
 
 \draw[thick] (2.5,2.5)--(2.5,3);
 \draw[thick] (2.5,-2.5)--(2.5,-3);

\node at (-1,0) {$\scalebox{1}{$k_1=$}$};

\end{tikzpicture}\qquad
\begin{tikzpicture}[x=.6cm, y=.6cm,
    every edge/.style={
        draw,
      postaction={decorate,
                    decoration={markings}
                   }
        }
]
\node at (-.5,0.5) {$\scalebox{.6}{$0$}$};
\node at (1,0.5) {$\scalebox{.6}{$2$}$};
\node at (3,0.5) {$\scalebox{.6}{$1$}$};
\node at (5,0.5) {$\scalebox{.6}{$1$}$};
\node at (4,0.5) {$\scalebox{.6}{$2$}$};
\node at (2,0.5) {$\scalebox{.6}{$0$}$};
\node at (1.5,0.2) {$\scalebox{.6}{$1$}$};
 
\node at (-.5,-0.5) {$\scalebox{.6}{$0$}$};
\node at (1,-0.5) {$\scalebox{.6}{$2$}$};
\node at (3,-0.5) {$\scalebox{.6}{$0$}$};
\node at (5,-0.5) {$\scalebox{.6}{$1$}$};
\node at (4.5,-0.25) {$\scalebox{.6}{$2$}$};
\node at (2.5,-0.2) {$\scalebox{.6}{$2$}$};
\node at (1.5,-0.25) {$\scalebox{.6}{$1$}$};
 \draw[dashed] (-1,0)--(6,0);
\node at (6,0.5) {\;};

 \draw[thick] (0,0)--(2.5, 2.5)--(5,0)--(2.5, -2.5)--(0,0) -- (1.5, 1.5)--(3,0)--(2,-1)--(1,0)--(2,1);
 \draw[thick] (1.5, .5) -- (2.5, -.5);
 \draw[thick] (2,-1) -- (3, -2);
 \draw[thick] (2,2) -- (4.5, -.5);
 
 \draw[thick] (2.5,2.5)--(2.5,3);
 \draw[thick] (2.5,-2.5)--(2.5,-3);

 \draw[thick] (6,-3)--(-1,-3);
 \draw[thick] (6,3)--(-1,3);

\end{tikzpicture}
\end{equation}

As for $M_2$, it is generated by $\CF$ and $\sigma(x_1)^2$. The strip depicted below shows that $\sigma(x_1)^2 \in \stab(S_1)$.
\[
\begin{tikzpicture}[x=.35cm, y=.35cm,
    every edge/.style={
        draw,
      postaction={decorate,
                    decoration={markings}
                   }
        }
]
 \node at (-11.5,0) {$\scalebox{1}{$\sigma(x_1)^2=x_0x_1x_2 x_0^{-3}=$}$}; 
 \node at (-6,0.5) {\;};

\draw[thick] (0,0) -- (-3,3)--(-6,0)--(-3,-3)--(0,0);
\draw[thick] (-1,1) -- (-4,-2);
\draw[thick] (-2,2) -- (-5,-1);

\draw[thick] (-3,3)--(-2,4)--(2,0)--(-2, -4)--(-3,-3);

\draw[thick] (-2,4)--(-2,4.5);
\draw[thick] (-2,-4)--(-2,-4.5);

\end{tikzpicture}
\qquad
\begin{tikzpicture}[x=.35cm, y=.35cm,
    every edge/.style={
        draw,
      postaction={decorate,
                    decoration={markings}
                   }
        }
]
\node at (2.5,0.5) {$\scalebox{.75}{$1$}$};
\node at (.5,0.5) {$\scalebox{.75}{$2$}$};
\node at (-1,0.5) {$\scalebox{.75}{$1$}$};
\node at (-3,0.5) {$\scalebox{.75}{$0$}$};
\node at (-5,0.5) {$\scalebox{.75}{$1$}$};
\node at (-7,0.5) {$\scalebox{.75}{$0$}$};

\node at (2.5,-0.5) {$\scalebox{.75}{$1$}$};
\node at (.5,-0.5) {$\scalebox{.75}{$2$}$};
\node at (-1,-0.5) {$\scalebox{.75}{$1$}$};
\node at (-3,-0.5) {$\scalebox{.75}{$2$}$};
\node at (-5,-0.5) {$\scalebox{.75}{$1$}$};
\node at (-7,-0.5) {$\scalebox{.75}{$0$}$};
\draw[dashed] (-7,0)--(3,0);
\node at (-6,0.5) {\;};

\draw[thick] (0,0) -- (-3,3)--(-6,0)--(-3,-3)--(0,0);
\draw[thick] (-1,1) -- (-4,-2);
\draw[thick] (-2,2) -- (-5,-1);

\draw[thick] (-3,3)--(-2,4)--(2,0)--(-2, -4)--(-3,-3);

\draw[thick] (-2,4)--(-2,4.5);
\draw[thick] (-2,-4)--(-2,-4.5);

\draw[thick] (-7,-4.5)--(3,-4.5);
\draw[thick] (-7,4.5)--(3,4.5);

\end{tikzpicture}
\]
By Remark \ref{color_at_extremes}, the color of the region to the right of the leftmost leaf is $2$ when the path from the root of the top/bottom tree to the leaf has odd length, and $1$ otherwise. 
For any $g$ in $\stab(S_1)$, these paths must have the same parity and thus $\log_2 g'(0)\in 2\IZ$. 
The subgroup $\stab(S_1)$ is bigger than $K_{(2, 2)}$. In fact, $\sigma(x_0x_1)\in \stab(S_1)\cap (K_{(2, 1)}\setminus K_{(2,2)})$
\[
\begin{tikzpicture}[x=.35cm, y=.35cm,
    every edge/.style={
        draw,
      postaction={decorate,
                    decoration={markings}
                   }
        }
]
\node at (-6.5,0) {$\scalebox{1}{$x_1x_0^{-2}=$}$}; 
 \node at (-6,0.5) {\;};

\draw[thick] (-2,2)--(-1,3)--(2,0)--(-1,-3)--(-2,-2);

\draw[thick] (0,0) -- (-1,1); 
\draw[thick] (0,2) -- (-1,1); 
\draw[thick] (-2,2)--(-4,0)--(-2,-2)--(0,0);
 \draw[thick] (-1,1) -- (-2,0)--(-3,-1);

 \draw[thick] (-1,3)--(-1,3.5);
 \draw[thick] (-1,-3)--(-1,-3.5);

\end{tikzpicture}\qquad
\begin{tikzpicture}[x=.35cm, y=.35cm,
    every edge/.style={
        draw,
      postaction={decorate,
                    decoration={markings}
                   }
        }
]
\node at (2.5,0.5) {$\scalebox{.75}{$1$}$};
\node at (.5,0.5) {$\scalebox{.75}{$0$}$};
\node at (-1,0.5) {$\scalebox{.75}{$1$}$};
\node at (-3,0.5) {$\scalebox{.75}{$2$}$};
\node at (-5,0.5) {$\scalebox{.75}{$0$}$};
 
\node at (2.5,-0.5) {$\scalebox{.75}{$1$}$};
\node at (.5,-0.5) {$\scalebox{.75}{$2$}$};
\node at (-1,-0.5) {$\scalebox{.75}{$1$}$};
\node at (-3,-0.5) {$\scalebox{.75}{$2$}$};
\node at (-5,-0.5) {$\scalebox{.75}{$0$}$};
 \draw[dashed] (3,0)--(-5,0);
\node at (-6,0.5) {\;};

\draw[thick] (-2,2)--(-1,3)--(2,0)--(-1,-3)--(-2,-2);

\draw[thick] (0,0) -- (-1,1); 
\draw[thick] (0,2) -- (-1,1); 
\draw[thick] (-2,2)--(-4,0)--(-2,-2)--(0,0);
 \draw[thick] (-1,1) -- (-2,0)--(-3,-1);

 \draw[thick] (-1,3)--(-1,3.5);
 \draw[thick] (-1,-3)--(-1,-3.5);

 \draw[thick] (-5,-3.5)--(3,-3.5);
 \draw[thick] (-5,3.5)--(3,3.5);

\end{tikzpicture}
\]
So far we have proved that $M_2\leq \stab(S_1) \cap K_{(2,2)}\leq K_{(2,2)}$. 
Thanks to the maximality of   $M_2$  in $K_{(2,2)}$, \cite[Theorem 4.6]{TV2},
  we only need to show that 
 $\stab(S_1) \cap K_{(2,2)}$ is a proper subgroup of $K_{(2,2)}$. We have already seen in \eqref{primoese} that
 $x_0^2$ is in 
  $K_{(2,2)}$
  but not in
   $\stab(S_1)$.
\end{proof}

\section{The $\vec{F}$-index}\label{sec_knots}
Recall that every oriented knot can be realized as $\vec{\CL}(g)$ for some $g\in \vec{F}$, \cite{A}. Therefore, one may consider the $\vec{F}$-index of a knot $K$, namely the smallest number of leaves required by each binary tree in a binary tree diagram $(T_+,T_-)$ of $\vec{F}$ such that $K$ is realized as $\vec{\CL}(T_+,T_-)$ (with the same orientation). 
For any oriented knot $K$, we denote by 
$K^r$ the knot with opposite orientation.

Planar graphs $\Gamma(T_+,T_-)$ arising from elements of $\vec{F}$ are 
rather special Tait graphs of links (see 
\cite[Proposition 4.1.3]{Jo14} or the discussion below the sole Theorem of \cite{A}).
We briefly recall from \cite{Jo14} that, given a link diagram $D$, its Tait graph is obtained as follows. First color the regions of the plane in black and white, 
with the unbounded region being black. The vertices of $\Gamma(D)$ sit in the interior of the black regions. 
Then we draw 
an edge between two vertices if and only if the two regions are the black regions of the same crossing. Each edge has a sign according to the following convention.
\[
\begin{tikzpicture}[x=.35cm, y=.35cm,
    every edge/.style={
        draw,
      postaction={decorate,
                    decoration={markings}
                   }
        }
]
\begin{scope}
\draw[thick] (2,0)--(0,2); 
\draw[line width = 3pt, white] (0,0) --(2,2);
\draw[thick] (0,0) --(2,2);

\fill[gray!20] (0,.1) -- (.9,1) -- (0,1.9) -- cycle;

\fill[gray!20] (2,.1) -- (1.1,1) -- (2,1.9) -- cycle;

\node at (-2,1) {$+=$};
\end{scope} 
\begin{scope}[xshift=3cm]
\draw[thick] (0,0) --(2,2);
\draw[line width = 3pt, white] (2,0)--(0,2);
\draw[thick] (2,0)--(0,2);

\fill[gray!20] (0,.1) -- (.9,1) -- (0,1.9) -- cycle;

\fill[gray!20] (2,.1) -- (1.1,1) -- (2,1.9) -- cycle;

\node at (-2,1) {$-=$};
\end{scope} 
\end{tikzpicture}
\]

Tait graphs completely determine a link and Reidemeister moves can be translated into this formalism.
The moves relevant to this paper read as shown in Figure \ref{ReideTait}.

\begin{figure}
\[
\begin{tikzpicture}[x=.35cm, y=.35cm,
    every edge/.style={
        draw,
      postaction={decorate,
                    decoration={markings}
                   }
        }
]
 
\begin{scope}
 \draw[thick] (-1,-1)--(0,0) -- (0,1) -- (-1,2) -- (-1.5,0.5)-- cycle;

 \draw[thick] (-1.5,0.5)--(-3,0.5);

\draw[fill=black] (0,0) circle (1pt); 
\draw[fill=black] (-1,-1) circle (1pt); 
\draw[fill=black] (-1,2) circle (1pt); 
\draw[fill=black] (0,1) circle (1pt); 
\draw[fill=black] (-1.5,0.5) circle (1pt); 
\draw[fill=black] (-3,0.5) circle (1pt); 
 
\node at (2.5,0.75)  {$\scalebox{.75}{$\stackrel{R1}{\longleftrightarrow}$}$};
\end{scope}
\begin{scope}[xshift=2.5cm]

 \draw[thick] (-1,-1)--(0,0) -- (0,1) -- (-1,2) -- (-1.5,0.5)-- cycle;

\draw[fill=black] (0,0) circle (1pt); 
\draw[fill=black] (-1,-1) circle (1pt); 
\draw[fill=black] (-1,2) circle (1pt); 
\draw[fill=black] (0,1) circle (1pt); 
\draw[fill=black] (-1.5,0.5) circle (1pt); 
  
 \end{scope}
\end{tikzpicture}
\]
\[
 \quad
\begin{tikzpicture}[x=.35cm, y=.35cm,
    every edge/.style={
        draw,
      postaction={decorate,
                    decoration={markings}
                   }
        }
]
\begin{scope}[xshift=1.5]
 \draw[thick] (-1,1)--(0,0);
 \draw[thick] (-1,0)--(0,0);
 \draw[thick] (-1,-1)--(0,0);
 
 \draw[thick] (0,0) to[out=-90, in=-90] (2,0);  
 \draw[thick] (0,0) to[out=90, in=90] (2,0);

 \draw[thick] (3,1)--(2,0);
 \draw[thick] (3,-1)--(2,0);

\draw[fill=black] (0,0) circle (1pt); 
\draw[fill=black] (2,0) circle (1pt); 

\node at (1,1.25)  {$\scalebox{.75}{$\pm$}$};
\node at (1,-1.25)  {$\scalebox{.75}{$\mp$}$};
\node at (5.25,0.25)  {$\scalebox{.75}{$\stackrel{R2}{\longleftrightarrow}$}$};
\end{scope}
\begin{scope}[xshift=3cm]
 \draw[thick] (-1,1)--(0,0);
 \draw[thick] (-1,0)--(0,0);
 \draw[thick] (-1,-1)--(0,0);

 \draw[thick] (3,1)--(2,0);
 \draw[thick] (3,-1)--(2,0);

\draw[fill=black] (0,0) circle (1pt); 
\draw[fill=black] (2,0) circle (1pt); 
 
\end{scope}
\end{tikzpicture}
\]
\[
\begin{tikzpicture}[x=.35cm, y=.35cm,
    every edge/.style={
        draw,
      postaction={decorate,
                    decoration={markings}
                   }
        }
]
 
\begin{scope}
 \draw[thick] (-1,1)--(0,0);
 \draw[thick] (-1,0)--(0,0);
 \draw[thick] (-1,-1)--(0,0);
 
 \draw[thick] (0,0) -- (1,0);  
 \draw[thick] (2,0) -- (1,0);

 \draw[thick] (3,1)--(2,0);
 \draw[thick] (3,-1)--(2,0);

\draw[fill=black] (0,0) circle (1pt); 
\draw[fill=black] (1,0) circle (1pt); 
\draw[fill=black] (2,0) circle (1pt); 

\node at (.5,.5)  {$\scalebox{.75}{$\pm$}$};
\node at (1.5,.5)  {$\scalebox{.75}{$\mp$}$};
\node at (5.25,0)  {$\scalebox{.75}{$\stackrel{R2}{\longleftrightarrow}$}$};
\end{scope}
\begin{scope}[xshift=3cm]
 \draw[thick] (-1,1)--(0,0);
 \draw[thick] (-1,0)--(0,0);
 \draw[thick] (-1,-1)--(0,0);
  
 \draw[thick] (1,1)--(0,0);
 \draw[thick] (1,-1)--(0,0);

\draw[fill=black] (0,0) circle (1pt); 
 
\node at (1,1.25)  {$\scalebox{.75}{$\pm$}$};
\node at (1,-1.25)  {$\scalebox{.75}{$\mp$}$};
\end{scope}
\end{tikzpicture}
\]
\caption{
The Reidemeister moves of type 1 and 2 in the language of Tait graphs.
In the second and third rows, $\pm$ and $\mp$ should be interpreted as meaning that if one edge has a positive sign, the other has a negative one. All unspecified signs are arbitrary.
}\label{ReideTait}
\end{figure}

More precisely a Tait graph is the planar graph of a certain element of $\vec{F}$ 
if and only if
 it satisfies the following properties (possibly after applying a planar isotopy) 
\cite[Proposition 4.1.3]{Jo14}  
\begin{enumerate}
\item it is bipartite;
\item the vertices are $(0,0)$, \ldots , $(N,0)$; 
\item each vertex is connected to exactly one vertex to its left;
\item each edge $e$ can be parametrized by a function $(x_e(t),y_e(t))$ with $x'_e(t)>0$, for all $t\in [0,1]$, and either $y_e(t)>0$, for all $t\in ]0,1[$ or $y_e(t)<0$, for all $t\in ]0,1[$;
\end{enumerate}

We conclude this paper with a result concerning the relation between the $\vec{F}$-index of a knot and its orientation.
\setcounter{theorem}{1}
\begin{theorem}\label{indextheo}
For any pair $(K, K^r)$ consisting of a knot and its reversed knot, the $\vec{F}$-indices of $K$ and $K^r$ differ by at most $4$.
\end{theorem}
\begin{proof}
Let $g \in \vec{F}$ be an element such that $\vec{\CL}(g) = K$, attaining the $\vec{F}$-index. We aim to construct an element $g'$ such that $\vec{\CL}(g') = K^r$.
There are two cases to consider: either the first and last vertices of the planar graph $\Gamma(g)$ associated with $g$ have the same color, or they have different colors.

In the first case, we insert the planar graph $\Gamma(g)$ as a subgraph of another graph
which still satisfies the properties characterizing this type of graphs:

\begin{center}
\begin{tikzpicture}[x=0.8cm,y=0.8cm,
  orig/.style={circle,draw=black,fill=red,inner sep=1.5pt},
  new/.style={circle,draw=black,fill=black,inner sep=1.5pt},
  lab/.style={below=7pt,font=\footnotesize}]
  \node[orig] at (0,0) {};  \node[lab] at (0,0) {$+$};
  \node at (1,0) {$\cdots$};
  \node[orig] at (2,0) {};  \node[lab] at (2,0) {$+$};
  \node at (3,0) {$\mapsto$};
  \def\A{4}\def\Vone{5}\def\VN{7}\def\Rone{8}\def\Rtwo{9}\def\Rthree{10}
  \draw (\A,0) .. controls (\A+0.25,0.575) and (\Vone-0.25,0.575) .. (\Vone,0);
  \draw (\VN,0) .. controls (\VN+0.25,0.575) and (\Rone-0.25,0.575) .. (\Rone,0);
  \draw (\A,0) .. controls (\A+1.25,2.4) and (\Rtwo-1.25,2.4) .. (\Rtwo,0);
  \draw (\Rtwo,0) .. controls (\Rtwo+0.25,0.575) and (\Rthree-0.25,0.575) .. (\Rthree,0);
  \draw (\A,0) .. controls (\A+0.25,-0.575) and (\Vone-0.25,-0.575) .. (\Vone,0);
  \draw (\VN,0) .. controls (\VN+0.25,-0.575) and (\Rone-0.25,-0.575) .. (\Rone,0);
  \draw (\Rone,0) .. controls (\Rone+0.25,-0.575) and (\Rtwo-0.25,-0.575) .. (\Rtwo,0);
  \draw (\VN,0) .. controls (\VN+0.75,-1.6) and (\Rthree-0.75,-1.6) .. (\Rthree,0);
  \node[new]  at (\A,0) {};      \node[lab] at (\A,0) {$+$};
  \node[orig] at (\Vone,0) {};   \node[lab] at (\Vone,0) {$-$};
  \node at (6,0) {$\cdots$};
  \node[orig] at (\VN,0) {};     \node[lab] at (\VN,0) {$-$};
  \node[new]  at (\Rone,0) {};   \node[lab] at (\Rone,0) {$+$};
  \node[new]  at (\Rtwo,0) {};   \node[lab] at (\Rtwo,0) {$-$};
  \node[new]  at (\Rthree,0) {}; \node[lab] at (\Rthree,0) {$+$};
\end{tikzpicture}
\end{center}

In this move, the first and last vertices in the initial graph (here depicted in red) change color
to $-$ and become the second and fourth to the last vertices of the new graph on the right. All the
other vertices and edges present in the initial graph are left unchanged (these are not drawn in the
figure). Overall, we added four vertices and eight edges (all depicted in black), which can all be
removed by means of Reidemeister moves of types I and II: two moves of type II remove the pair of
parallel edges joining the first two vertices and the pair of parallel edges joining the last vertex
of the initial graph to the following one, after which the four added vertices become, one after the
other, univalent and are removed by moves of type I. This means that topologically the knots are the
same; however, the colors of the regions are the opposite and so is the orientation of the knot.

When the first and last vertices have opposite colors,  it suffices to apply the following move. 
As in the other move, the first and the last vertices in the initial graph (appearing in red) both change color and become the second and second to the last vertices of the new graph on the right. All the other vertices and edges present in the initial graph are left unchanged.
In total, we added two vertices and four edges (all depicted in black), which can all be removed by means of Reidemeister moves of types I and II. This means that topologically the knots are the same; however, the colors of the regions are the opposite and so is the orientation of the knot.
\[
\begin{tikzpicture}

        \draw[fill=red] (0,0) circle (2pt); 
        \draw[fill=red] (2,0) circle (2pt); 
        \node at (1,0) {$\scalebox{.75}{$\ldots$}$};
        \node at (3,0) {$\scalebox{.75}{$\mapsto$}$};
        \node at (0,-.25) {$\scalebox{.75}{$+$}$};
        \node at (2,-.25) {$\scalebox{.75}{$-$}$};
\end{tikzpicture}
\quad
\begin{tikzpicture}
    \draw (-1,0) to[out=90, in=90] (3,0);  
    \draw (-1,0) to[out=90, in=90] (0,0);  
    \draw (-1,0) to[out=-90, in=-90] (0,0);  
    \draw (2,0) to[out=-90, in=-90] (3,0);

        \draw[fill=black] (-1,0) circle (2pt); 
        \draw[fill=red] (0,0) circle (2pt); 
        \draw[fill=red] (2,0) circle (2pt); 
        \draw[fill=black] (3,0) circle (2pt); 
         \node at (1,0) {$\scalebox{.75}{$\ldots$}$};
         \node at (0,-.25) {$\scalebox{.75}{$-$}$};
        \node at (2,-.25) {$\scalebox{.75}{$+$}$};
        
         \node at (-1,-.25) {$\scalebox{.75}{$+$}$};
        \node at (3,-.25) {$\scalebox{.75}{$-$}$};
\end{tikzpicture}
\]
These moves do not affect  the knot topologically (they are in fact type I and type II Reidemeister moves), they only change its orientation.
\end{proof}
We end the paper with the following natural questions.
\begin{question}
Does there exist a pair $(K, K^r)$ consisting of an oriented knot and its reverse that are equivalent as unoriented knots but have different $\vec{F}$-indices?
\end{question}
\begin{question}
In the proof of Theorem \ref{indextheo} we have seen that the difference between the $\vec{F}$-indices of $K$ and $K^r$ can be at most $2$ when the minimal planar graph has a coloring where the first and last vertices have different colors, while the difference is at most $4$ when the first and last vertices have the same color.
Does every knot have minimal diagrams (with the same number of vertices) of both types? If not, do the knots of the two types differ from each other?
\end{question}

\section*{Acknowledgements} 
The motivation for writing this paper stems from the questions asked during two presentations: one given in 2024 at the Operator Algebras seminar at the University of Rome 'Tor Vergata', during which Tiziano Gaudio posed the question of whether \( M_0 \), \( M_1 \), and \( M_2 \) could be stabilizers of certain subsets of dyadic rationals; the other given in 2021 at the Geometric Topology seminar of Columbia University on Jones's construction of knots. We extend our gratitude to T.G. for his interest in this topic. Finally, I would like to thank the referee for their perusal of the manuscript and for comments that resulted in changes making the article more accessible and self-contained. The author is partially supported by Sapienza Università di Roma (Progetto di Ateneo Dipartimentale 2024 “New research trends in Mathematics at Castelnuovo”) and by INdAM-GNAMPA through Project "Simmetrie distribuzionali per processi stocastici quantistici" CUP E53C25002010001.
   
%

%


\end{document}